\newif\ifmakeplots
\makeplotstrue

\newif\ifarxiv
\arxivtrue

\newif\ifjota
\jotafalse

\ifjota
    \documentclass[smallextended,envcountsect]{svjour3} 
\else
    \documentclass[12pt,notitlepage]{article}
\fi

\ifarxiv
    \newcommand{\datadir}{.}
    \newcommand{\ext}{.txt}
\else
    \newcommand{\datadir}{data}
    \newcommand{\ext}{}
\fi

\usepackage[anythingbreaks]{breakurl}

\ifjota
\else
    \usepackage[margin=1.25in, headheight=42pt]{geometry}
    \usepackage{amsthm}
\fi
\usepackage{amsmath}
\usepackage{thmtools}
\usepackage{booktabs}
\usepackage{cleveref}
\usepackage[hyphens]{url}
\usepackage{graphicx}
\usepackage{pgfplots}
\usepackage{algorithm2e}
\usepackage{environ}
\usepackage{dcolumn}
\usepackage{filecontents}
\ifjota
\else
    \usepackage{subcaption}
\fi
\usepackage[font=small,labelfont=bf]{caption}
\usepackage[style=authoryear,backend=bibtex,giveninits=true,url=false,doi=false,dashed=false,natbib=true,maxbibnames=99]{biblatex}
\usepackage{xurl}

\newcolumntype{d}{D{.}{.}{-1}}

\usepackage{fancyhdr}

\ifjota
\else
    \fancypagestyle{plain}{%
        \fancyhead{}
        \chead{\textsc{For professional clients and qualified investors only\\Not for public distribution\\This version posted with permission}}
    }
\fi

\usepgfplotslibrary{groupplots}
\pgfplotsset{
    layers/my layer set/.define layer set={
        background,
        main,
        foreground
    }{
    },
    set layers=my layer set,
}

\definecolor{aux_color}{rgb}{1, 0, 0}  
\definecolor{soph_heur_color}{rgb}{0, 0, 1}  
\definecolor{basic_heur_color}{rgb}{0, .5, 0}  
\definecolor{boxcolor}{rgb}{0.5, 0.5, 0.5}  

\def\auxcolorname{red}
\def\sophcolorname{blue}
\def\basiccolorname{green}

\usetikzlibrary{
    pgfplots.dateplot,
}
 
\tikzset{
    font={\small\selectfont}
}

\pgfplotsset{compat=newest}

\def\aspectratio{1.3}

\newcommand{\eg}{{\it e.g.}}
\newcommand{\ie}{{\it i.e.}}

\newcommand{\BA}{\begin{array}}
\newcommand{\EA}{\end{array}}

\newcommand{\ones}{\mathbf 1}

\newcommand{\reals}{{\mbox{\bf R}}}

\newcommand{\sign}{\mathop{\bf sign}}

\title{Tax-Aware Portfolio Construction\\via Convex Optimization}

\ifjota
    \author{Nicholas Moehle \and Mykel J. Kochenderfer \and Stephen Boyd \and Andrew Ang}
\else
    \author{Nicholas Moehle \and Mykel J. Kochenderfer \and Stephen Boyd \and Andrew Ang\\[10pt]
    BlackRock AI Labs}
\fi

\begin{document}
	
\maketitle

\begin{abstract}
We describe an optimization-based tax-aware portfolio construction method that 
adds tax liability to standard Markowitz-based portfolio construction.
Our method produces a trade list that specifies
the number of shares to buy of each asset and
the number of shares to sell from each tax lot held.
To avoid wash sales (in which some realized capital losses are disallowed),
we assume that we trade monthly, and cannot simultaneously buy and sell the
same asset.

The tax-aware portfolio construction problem is not convex,
but it becomes convex when we specify, for each asset, whether we buy or sell it.
It can be solved using standard mixed-integer convex optimization methods
at the cost of very long solve times for some problem instances.
We present a custom convex relaxation of the problem
that borrows curvature from the risk model. This relaxation can provide a good
approximation of the true tax liability, while greatly enhancing computational tractability.
This method requires the solution of 
only two convex optimization problems: the first determines
whether we buy or sell each asset, and the second generates the 
final trade list.
In our numerical experiments, 
our method almost always solves the nonconvex problem to optimality,
and when it does not, it produces a trade list very close to optimal.
Backtests show that the performance of our method is indistinguishable 
from that obtained using
a globally optimal solution, but with significantly reduced computational effort.
\end{abstract}

\clearpage

\section{Introduction}

We formulate a tax-aware portfolio construction problem
that explicitly accounts for tax liabilities from long- and short-term capital gains, while
maintaining the standard objectives of return, risk, and transaction costs.
While this tax liability term is not convex,
we develop a convex relaxation of the problem that allows 
trade lists to be computed extremely quickly.

In addition to choosing how many shares to buy or sell of each asset,
our method divides each sale across tax lots.
For a fixed sell amount,
we find it is optimal to sell shares in the order that minimizes immediate tax liability,
which we call the least tax first out (LTFO) method,
coinciding with the well-known highest basis first out (HIFO) method if
the same tax rate applies to all lots.
\citet{Dickson2000} and \citet{Berkin2003} show that the accounting method for lot
ordering significantly affects the losses that can be harvested.
Even though it is not optimal,
\citet{Atra2014} show that the HIFO method generates
substantial benefit to an investor's total wealth.

Our convex relaxation combines all terms that are separable across a single asset,
and then replaces the resulting (nonconvex) function with its convex envelope.
This approach has the effect of taking (or `borrowing') curvature
from the transaction cost and specific risk terms
and introducing it into the tax liability term.
This effect can also be interpreted as an application of the Shapley--Folkman Lemma.
The effect is that the tax-aware portfolio construction problem can
be well approximated by a convex optimization problem,
and the fidelity of this approximation is shown empirically.


To demonstrate our method, we apply it to a tax-loss harvesting strategy.
First, we show that in a realistic backtest scenario,
the strategy tightly tracks the benchmark (the S\&P 500)
while harvesting capital losses.
We compare the performance of our method with that of a standard mixed-integer quadratic programming formulation
and show that our method delivers near-identical trade lists,
despite being several hundred times faster.

\subsection{Related work and background}

\paragraph{Markowitz portfolio construction.}
The formulation of portfolio construction
as an optimization problem by \citet{markowitz1952} involves a trade
off of expected return and risk.  This optimization problem,
with a quadratic objective and linear equality constraints,
has an analytical solution. 
The problem can be extended by including position limits or a long-only constraint 
\citep{sharpe1963simplified,
markowitz1955optimization,
grinold2000active}. 
The resulting problem no longer has an analytical solution,
but it can be efficiently solved as a quadratic program (QP) \citep[pp.~55--156]{cvxbook}.
Other constraints and objective terms can be incorporated, including those related to accounting for a previous or 
initial portfolio \citep{pogue1970extension, lobo2007portfolio}.
Including an initial portfolio allows the method to
be used as a trading policy,
which can be run periodically to prescribe trades \citep{boyd2017multi}.

\paragraph{Portfolio construction via convex optimization.}
These Markowitz-inspired portfolio construction problems are typically convex
and can be efficiently solved \citep{cvxbook}.
Even complex portfolio construction problems can be specified succinctly
in a high-level domain-specific language for convex optimization,
such as CVXPY \citep{diamond2016cvxpy}, CVX \citep{cvx}, Convex.jl \citep{convexjl}, and 
CVXR \citep{cvxr}.
Furthermore, problems with thousands of assets and a risk model with dozens of 
factors can be solved in well under a second using standard open-source solvers
such as ECOS \citep{domahidi2013ecos},
OSQP \citep{osqp}, SCS \citep{ocpb:16, scs}, or commercial solvers 
such as CPLEX \citep{cplex}, MOSEK \citep{mosek}, or GUROBI \citep{gurobi}.
Custom implementations of portfolio construction solvers can be far faster, with 
solve times measured in milliseconds.

Although solver speed is not essential if we are only interested in occasionally rebalancing a handful of portfolios,
it can be useful when managing a very large number of individualized accounts.
Additionally, fast solvers allow us to quickly run many backtest simulations of a trading algorithm.
These simulations allow us to tune hyper-parameters, carry out what-if experiments, and compare
different formulations or models on historical or synthesized data.

\paragraph{Non-convex portfolio construction problems.}
Some practical constraints and objective terms are not convex.
An obvious example is that asset holdings must be 
in integral numbers of shares.
For large portfolios, this constraint is readily 
handled by simple heuristics, for example by ignoring it in
solving the problem, and then rounding the real-valued
holdings to the nearest integer values.
Other more challenging constraints include limits on the number of 
assets in the portfolio, or a minimum nonzero trade size
\citep{bertsimas1999portfolio}.
A challenging nonconvex objective term is tax liability, the
focus of this paper.

These nonconvex portfolio construction problems can be reformulated as 
mixed-integer convex optimization problems (which are not convex).
They can be solved exactly using a variety of methods and software,
such as GLPK \citep{glpk}, CPLEX, MOSEK, and GUROBI.
Such solvers are often fast, but for 
some problem instances can have very long solve times,
often hundreds of times more than those associated with similar convex problems.
In contrast, solving convex optimization problems is reliably fast.  

\paragraph{Convex approximations.} 
An alternative to solving the nonconvex optimization problem exactly is 
to employ a heuristic method that finds an approximate solution far 
faster than it would take to solve the problem exactly.
This paper presents one such heuristic method, based on
a convex approximation of the original problem.
The idea that convex approximations of nonconvex problems can be used
in place of global nonconvex solvers with the same practical performance has been widely noted in other areas
\citep{NCVX}. 

The problem we study in this paper
involves the sum of many nonconvex terms that are all similar.
Because the sum of a large number of nonconvex functions tends to be `more convex'
than the original functions,
these problems are often well approximated by convex problems.
This intuitive phenomenon was described by Shapley, Folkman, and Starr.
Starr applied it to problems in microeconomics involving many agents \citep{starr1969quasi}.
\citet{bertsekas1997nonlinear} provides practical algorithms for solving these problems with performance bounds.

\paragraph{Tax-aware investment.}
Our paper is related to literature that develops optimal tax-aware
trading strategies,
building on the papers by \citeauthor{Constantinides1983} \citeyearpar{Constantinides1983,Constantinides1984}.
The intuition in these papers
is that investors may reduce their tax liability by deferring capital gains and realizing losses,
which can be used to offset current income or capital gains.
This intuition also applies in our setting.
Many papers apply well-known numerical
methods to solve the tax problem.
For example, dynamic programming has been widely applied \citep{Dammon1996,Dammon2004,Dammon2001}.
\citet{Dybvig1996} use a
binomial tree and formulate an optimal stopping problem. 
\citet{DeMiguel2005} formulate an optimal tax investment strategy with nonlinear
programming.
Although Markowitz-based tax-aware portfolio construction is an old idea
\citep{pogue1970extension},
to our knowledge, ours is the first to develop a \emph{convex} tax optimization problem by relaxing
the original nonconvex problem.
We focus on the speed and
reliability of convex optimization techniques with applications to 
taxable managed funds \citep{Sialm2020}, the large and rapidly growing
tax-loss harvesting industry \citep{Chaudhuri2020}, and
security valuation with taxes \citep{Gallmeyer2011}.

\paragraph{Multi-period portfolio construction.}
We focus on the single-period portfolio construction problem,
without explicitly planning for future trades.
This is in contrast to multi-period portfolio optimization formulations,
such as those of \citet{boyd2017multi}.
Many tax-aware problems are readily handled by single-period portfolio optimization,
such as loss harvesting, tax-neutral portfolio rebalancing, managing inflows and outflows,
and optimizing tax-free donations.
In fact, we believe repeated single-period portfolio optimization is a excellent approach
to certain long-horizon investment problems,
such as tracking a (low-turnover) benchmark portfolio while harvesting tax losses.
There is some theoretical justification here:
\citet{Constantinides1983} shows that for a single asset,
with no wash sale rule, the greedy approach of realizing losses and deferring gains is optimal.
Repeatedly using single-period, tax-aware portfolio construction does exactly this,
but with many assets, while avoiding wash sales.
A more practical justification for single-period portfolio optimization 
is that it is currently standard for such tax-loss-harvesting strategies.
This paper does not argue for or against using a single-period formulation,
but simply gives a reasonable problem formulation and solution method for it.

\subsection{Contributions}

This paper focuses on incorporating a specific nonconvex term, the tax liability
generated by the trades, into an otherwise convex portfolio construction problem.
Ignoring this constraint, or using simple ad hoc
rounding methods to handle it, does not work well compared to solving
the problem exactly with a mixed-integer convex solver.
Our contribution is to develop a heuristic method for approximately solving
the tax-aware portfolio construction problem that relies on solving 
two convex optimization problems, making it reliably fast.

\section{Tax-aware portfolio construction}\label{s-TAM}

This section outlines our notation and describes the tax-aware portfolio
optimization problem.
We start by describing the trading dynamics and various objective terms.

\subsection{Portfolio holdings and dynamics}
We consider a universe of $n$ assets we are allowed to hold and trade.
We let $h_{\rm init} \in \reals^n$ denote the dollar value of our pre-trade holdings
of these $n$ assets.
We restrict ourselves to long-only portfolios, so $h_{\rm init} \geq 0$.

Our task is to decide how much of each these assets to buy or sell.
We represent this decision by 
a purchase vector $u\in\reals^n$, denominated in dollars.
If we purchase asset $i$, $u_i > 0$; if we sell asset $i$, $u_i < 0$.
Our post-trade holdings are $h\in\reals^n$, given by
\[
h = h_{\rm init} + u.
\]
This equation ignores transaction costs, which are assumed to be small.
(Following convention,
we include these transaction costs in our objective function.)
We require that the post-trade portfolio is also long-only, making $h \ge 0$.
This constraint means we cannot sell more of any asset than we 
currently hold.

\paragraph{Cash.}
The cash held in the portfolio is $c_{\rm init}\in\reals$,
which we allow to be negative.
The post-trade cash balance is
\[
c = c_{\rm init} - \ones^Tu.
\]
We assume the post-trade cash amount must match some desired value $c_{\rm des}$, which translates to the constraint on $u$
\[
\ones^T u = c_{\rm init}-c_{\rm des}.
\]

The total pre-trade portfolio value, including cash,
is $\ones^T h_{\rm init} + c_{\rm init}$, which we assume is positive.
While any value of $c_{\rm des}$ is possible, a common choice
is a given fraction $\eta$ of the total portfolio value, 
\begin{equation}\label{e-cdes-choice}
c_{\rm des}=\eta(\ones^T h_{\rm init} + c_{\rm init}).
\end{equation}
The choice $\eta=0.01$, for example, means that 1\% of the total portfolio 
value is to be held in cash.
The cash balance can be used to handle cash deposits into and withdrawals from 
the account by adjusting $c_{\rm init}$ by the amount deposited or 
withdrawn.

\subsection{Objective terms}
Here we describe various objective terms and additional constraints,
including the traditional ones: expected return, active risk, and 
transaction costs.  
We briefly introduce the tax liability term and provide some of its
attributes, reserving a detailed description for section~\ref{s-tax-liab}.

We note that it is customary to scale the variables and objective terms
so that the objective represents an adjusted return.
However, to simplify our description of the tax liability function in section~\ref{s-tax-liab},
we leave the variables and objective in units of dollars.

\paragraph{Risk.}
The risk of a managed portfolio is typically measured
with respect to a benchmark portfolio, such as the S\&P 500.
This benchmark portfolio is described by a vector $h_b\in\reals^n$,
scaled so that it has the same market value as our portfolio,
\ie, $\ones^T h_b = \ones^T h_{\rm init} + c_{\rm init}$.

The (active) risk is
\[
(h - h_b)^T V (h - h_b),
\]
where $V$ is the covariance matrix of the asset returns.
Our covariance matrix $V$ has the traditional factor model form
\[
V = X\Sigma X^T + D,
\]
where $X\in \reals^{n \times k}$ is the factor exposure matrix,
$\Sigma\in \reals^{k \times k}$ is the symmetric positive definite 
factor covariance
matrix, and $D\in\reals^{n\times n}$ is the diagonal matrix of idiosyncratic variances with $D_{ii} > 0$
\citep{grinold2000active, boyd2017multi}.
The risk can be decomposed into two components,
the systematic risk
\begin{equation}\label{e-factor-risk}
(h - h_{b})^T X\Sigma X^T (h - h_{b}),
\end{equation}
and the specific risk
\begin{equation}\label{e-idio-risk}
(h - h_{b})^T D (h - h_{b}) = \sum_{i=1}^n D_{ii} (h_i-h_{b,i})^2.
\end{equation}
It is common to express active risk in terms of its square root, which has 
units of dollars.
We note for future use that the specific risk~(\ref{e-idio-risk}) is 
separable, \ie, a sum of terms each associated with one asset.

\paragraph{Expected return.}
Suppose we have a forecast of the return of the $n$ assets,
expressed as a vector $\alpha\in\reals^n$,
where $\alpha_i$ is the expected return of asset $i$.
The expected active return
of portfolio $h$ is then $\alpha^T (h - h_b)$, which is measured in dollars.
Because a constant offset is immaterial for optimization,
we can write the expected return as simply $\alpha^T h$ or even $\alpha^T u$.

\paragraph{Transaction costs.}
The transaction cost follows a simple bid-ask spread model:
\[
\kappa^T |u|,
\]
where $\kappa\in\reals^n_{+}$ is one-half the bid-ask spread,
and $|u|$ is the element-wise absolute value of $u$.
For simplicity, we neglect the standard price impact term;
this omission is reasonable if we assume our trades
are small relative to the total market volume over the trading period.
(For larger accounts, a price impact term can be included; see \citet{boyd2017multi}.)

\paragraph{Tax liability.}
We let $L: \reals^n \to \reals$ denote the tax liability function,
where $L(u)$ is the immediate tax liability incurred by the trades $u$ due to realizing capital gains.
We will describe $L(u)$, which derives from the history
of previous transactions in the assets, in detail 
in section~\ref{s-tax-liab};
for now, we simply note some of its attributes.
First, it is separable across the assets, \ie, it has the form
\[
L(u) = \sum_{i=1}^n L_i(u_i),
\]
where $L_i(u_i)$ is the tax liability for asset $i$ incurred by trading.
There is no immediate tax liability when buying an asset, making $L_i(u_i)=0$
for $u_i\geq 0$.  For $u_i<0$, \ie, selling the asset, $L_i(u_i)$ is 
a convex piecewise linear function.
While $L_i$ is convex for $u_i<0$, it is not convex over the 
whole interval, which includes buying ($u_i>0$) and selling ($u_i<0$).
The total tax liability function $L(u)$ is not convex,
but it becomes convex if we restrict the sign of $u_i$, \ie, 
we specify whether we are buying or selling each asset.

\paragraph{Constraints.}
We have already mentioned several constraints, for example that $h \geq 0$
(the portfolio is long-only) and $\ones^T u = c_{\rm init}- c_{\rm des}$
(the post-trade cash matches a desired value).
We also allow for additional convex constraints on
the trade list and post-trade holdings. 
We represent these as $u \in \mathcal U$ and $h \in \mathcal H$.
These could include, for example, limits on the holdings of a particular asset,
or limits on the exposure of our portfolio to a certain factor.
For concreteness, we assume $\mathcal H$ and $\mathcal U$ are polyhedral,
\ie, described by a finite set of linear equality and inequality constraints,
although our proposed method also applies more generally.

\subsection{Tax-aware portfolio construction}

\paragraph{Tax-aware utility function.}
We assemble our objective terms into a single utility 
function of $u$ and $h$, 
\begin{align}
U(h,u)=
\alpha^T u
- \gamma_\text{risk} (h - h_{b})^T V (h - h_{b})
- \gamma_\text{tc} \kappa^T |u|
- \gamma_\text{tax} L(u).
\label{e-utility}
\end{align}
where $\gamma_\text{risk}$, $\gamma_\text{tc}$, and $\gamma_\text{tax}$
are nonnegative trade-off parameters.
The first two terms constitute the traditional risk-adjusted return used in Markowitz portfolio construction.
The third term is transaction cost, a widely used addition to the 
traditional Markowitz utility,
with the parameter $\gamma_\text{tc}$ used to control turnover.
The last term accounts for the tax liability of the 
trades.

\paragraph{Tax-aware portfolio construction problem.}
Our problem is to maximize utility subject to constraints, \ie,
\begin{equation}\label{e-prob}
\begin{array}{ll}
\text{maximize} &
\alpha^T u
- \gamma_\text{risk} (h - h_b)^T V (h - h_b)
- \gamma_\text{tc} \kappa^T |u|
- \gamma_\text{tax} L(u)  \\
\text{subject to} &
h = h_{\rm init} + u, \quad
\ones^T u = c_\text{init}-c_\text{des} \\
& u \in \mathcal U, \quad h \in \mathcal H,
\end{array}
\end{equation}
with decision variables $u$ and $h$.
The problem data are
$\alpha$,
$h_b$,
$V$,
$\kappa$,
$h_{\rm init}$,
$c_{\rm des}$,
$c_{\rm init}$,
the function $L$ (described in section~\ref{s-tax-liab}),
the constraint sets $\mathcal U$ and $\mathcal H$,
and the trade-off parameters
$\gamma_\text{risk}$,
$\gamma_\text{tc}$,
and
$\gamma_\text{tax}$.
We refer to the problem~(\ref{e-prob}) as the 
\emph{tax-aware Markowitz problem}, or TAM problem,
and we denote its optimal value as $U^\star$.

\paragraph{Non-convexity.}
The constraints in the TAM problem~(\ref{e-prob}) are convex,
as are all terms in the objective with the exception of the 
tax liability term.  Unfortunately, that term renders the
TAM problem (\ref{e-prob}) nonconvex,
which makes it difficult to solve (exactly) in general.
We note, however, that the problem becomes convex when we specify the 
sign of the trade list $u$,
\ie, if we specify for each asset whether we 
are to sell it ($u_i\leq 0$) or buy it ($u_i\geq 0$).

The TAM problem can be formulated as a mixed-integer quadratic program (MIQP), 
which can be solved using various methods.  It is well known that in practice,
these methods can often solve problems reasonably quickly, but in many 
other cases, the solution times can be extremely long.
The main contribution of this paper is a method for approximately
solving the TAM problem, which involves solving only two
convex optimization problems.  As a result, our method is always very fast
and never involves the very long solution time that can be observed with
MIQP solvers.
As we will see in section~\ref{s-examples},
for realistic instances of the TAM problem,
our method delivers near identical performance
as a globally optimal solution.

\section{Tax liability}\label{s-tax-liab}

This section describes the tax liability function.

\subsection{Tax lots and capital gains}

\paragraph{Tax lots.}
For each asset,
the pre-trade holdings are composed of zero or more tax lots.
Each tax lot has several attributes associated with it:
its quantity of shares, acquisition date,
and cost basis (the price per share at which the shares were acquired).

We let $q_{ij}$ denote the quantity of shares in the $j$th lot of asset $i$.
The total number of shares of asset $i$ held is $\sum_j q_{ij}$.
We have $h_{{\rm init}, i} = p_i \sum_j q_{ij}$,
where $p_i$ is the current price of asset $i$.
We let $b_{ij}$ denote the cost basis (in dollars per share)
of the $j$th lot of asset $i$.

\paragraph{Selling shares.}
When shares of asset $i$ are sold, \ie, we have $u_i<0$,
we must specify which tax lots from which to take the shares.
Let $s_{ij}$ denote the dollar value of shares sold from the $j$th lot of asset $i$,
with $0 \leq s_{ij} \leq q_{ij}p_i$, where $q_{ij}p_i$ is the dollar value
of the $j$th lot of asset $i$.
The total dollar value of asset $i$ sold is then $\sum_j s_{ij}$,
which must be equal to $-u_i$.

When we sell $s_{ij}$ dollars from lot $j$ of asset $i$,
we incur a capital gain, which is the difference of our proceeds
and our cost basis for those shares, \ie,
$(1 - b_{ij}/p_i)s_{ij}$.  We refer to this quantity as the gain;
when it is negative, we refer to it as the loss.

\paragraph{Long-term and short-term gains.}
A tax lot is \emph{long term} if the acquisition date is more than one
year before the trade date,
and the lot is \emph{short term} otherwise.
Gains from long-term and short-term lots
are taxed at two different positive rates,
$\rho_\text{lt}$ and $\rho_\text{st}$, respectively,
with $\rho_\text{lt} \le \rho_\text{st}$.
The tax liability for selling dollar value $s_{ij}$ from the $j$th lot of asset
$i$ is $\rho_\text{lt} (1- b_{ij}/p_i)s_{ij}$ if lot $j$ is long term,
and 
$\rho_\text{st}(1 - b_{ij}/p_i)s_{ij}$ if lot $j$ is short term.

The total tax liability from selling all assets is
\[
\sum_{i,j} \rho_{ij} (1 - b_{ij}/p_i) s_{ij} =
\sum_{i,j} T_{ij} s_{ij},
\]
where the tax rates $\rho_{ij}$ are given by: 
\[
\rho_{ij} =
\begin{cases}
\rho_\text{lt} & \text{lot $j$ of asset $i$ is long term} \\
\rho_\text{st} & \text{lot $j$ of asset $i$ is short term}\text.
\end{cases}
\]
We refer to $T_{ij} = \rho_{ij} (1- b_{ij}/p_i)$ as the
tax rate for lot $j$ of asset $i$.
This is the dollar tax liability generated per dollar sold of the lot.
It is positive if the current asset price exceeds the lot basis 
\ie, the lot is held at a gain,
and is negative if the lot is held at a loss.

\subsection{Tax liability function}
Suppose that for asset $i$ we have $u_i < 0$,
\ie, we are selling $-u_i$ dollars of asset $i$,
which translates to $-u_i / p_i$ shares.
We can solve the problem of allocating the sale across
lots in order to minimize the tax liability incurred.
We define
\[
L_i(u_i) = \min_{s_{ij}} \left\{ \sum_j T_{ij}s_{ij} \;\middle|\;
\sum_j s_{ij} = -u_i , ~0 \leq s_{ij} \leq q_{ij}p_i \right\},
\]
which is the smallest tax gain achievable to carry out this sale.
We define $L_i(u_i) = +\infty$ for $-u_i < p_i\sum_j q_{ij}$,
\ie, if we ask to sell more shares of the asset than we hold.
We also define $L_i(u_i) = 0$ for $u_i \geq 0$,
\ie, we are buying shares instead of selling.
These properties hold because purchasing additional shares incurs no immediate tax liability.

\paragraph{Least-tax-first-out lot policy.}
For a given value of $u_i < 0$,
it is easy to determine optimal values of $s_{ij}$;
it is a convex optimization problem with an analytical solution.
We simply sort the values of $T_{ij}$ from
least (most negative) to greatest, breaking ties arbitrarily.
Then we sell shares from lots in this order.
For example, we start by selling shares from
the lot with the smallest (or most negative) value of
$T_{ij}$, which is the term-adjusted tax liability rate.
If we need to sell more shares than that, we go to the lot with
second smallest value, and so on.
This greedy approach is optimal, \ie, it minimizes
the tax liability when selling $-u_i$ dollars of asset $i$.
We refer to this approach of choosing lots from
which to sell shares as \emph{least tax first out} (LTFO).
It takes into account whether the lots are long term or short term.
If all lots are short term, this scheme
corresponds to the well-known \emph{highest basis first out} (HIFO) method.

The tax liability function $L_i$ is continuous and piecewise affine.
If none of the lots are at a loss, $L_i$ is convex and nonnegative.
If at least one lot is held at a loss,
then $L_i$ takes on negative values and is \emph{not} convex.
When the domain of $L_i$ is
restricted to either $u_i \le 0$ or $u_i \ge 0$,
the resulting function is convex.

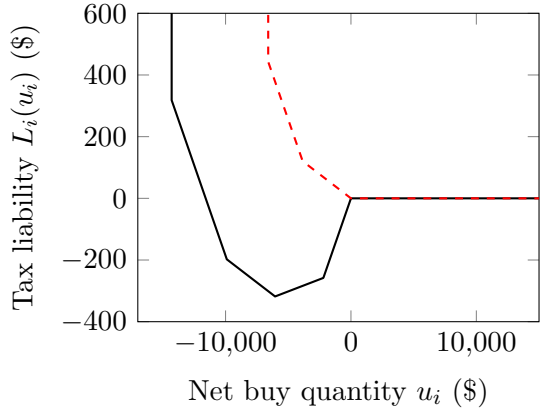
\begin{figure}
\centering
\ifmakeplots
    \begin{tikzpicture}
        \begin{axis}[xlabel = Net buy quantity $u_i$ (\$),
                     ylabel = Tax liability $L_i(u_i)$ (\$),
                     xmin=-17000,
                     xmax=15000,
                     ymin=-400,
                     ymax=600,
                     width = 0.35\textwidth,
                     height = 0.35\textwidth / \aspectratio,
                     scale only axis=true,
                     xticklabel style={/pgf/number format/fixed},
                     scaled ticks = false,
                     xtick={-10000, 0, 10000},
                     legend cell align={left}
                     ]

            \addplot [thick, black] table{\datadir/tax_liability_function_with_loss.csv\ext};

            \addplot [thick, dashed, aux_color] table{\datadir/tax_liability_function_with_gain.csv\ext};


        \end{axis}
    \end{tikzpicture}
\fi
\caption{Tax liability functions $L_i$ for two assets.
The solid black curve is for an asset with four lots
with two held at a loss and two held at a gain.
The dashed \auxcolorname{} curve is for an asset with two lots, both held at a gain.}
\label{f-taxliability-function}
\end{figure}

Figure~\ref{f-taxliability-function}
shows two different tax liability functions.
The dashed \auxcolorname{} curve shows the tax liability function for an asset
which we hold in two lots, both at a gain
(\ie, with current price greater than basis).
The solid black curve shows the tax liability function of a different asset,
which we hold in four lots, two at a loss
(\ie, the basis is greater than the current price),
and two at a gain.
Each linear segment corresponds to a tax lot,
with the slope given by the tax rate of the lot,
and width given by the total value of the lot.

\section{Convex relaxation}\label{s-conv-relaxations}

The first step in developing our convex-optimization-based 
heuristic for approximately solving the TAM problem is to
form a convex relaxation that approximates the problem.

\subsection{Convex envelope of a function}

We review a standard concept, the \emph{convex envelope} of a function $f:\reals
\to \reals$, denoted $f^{**}$.
It is defined as
\begin{equation}\label{e-conv-env-def}
f^{**}(x) = 
\inf \{\theta f(v) + (1-\theta) f(w) \mid \theta \in [0,1],~x = \theta v
+ (1 - \theta) w\}.
\end{equation}
The infimum is over $\theta$, $v$, and $w$.
The convex envelope function $f^{**}$ is convex,
and it satisfies $f^{**}(x) \leq f(x)$ for all $x$, \ie, it is a global
underestimator of $f$.
If $f$ is convex, then $f^{**}$ is equal to $f$.

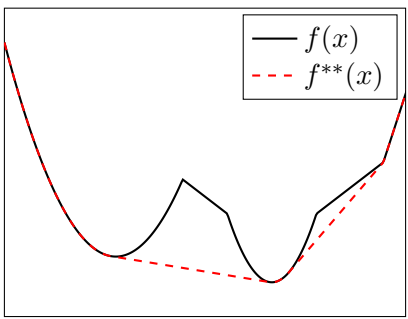
\begin{figure}
\centering
\ifmakeplots
    \begin{tikzpicture}
        \begin{axis}[xmin=-1,
                     xmax=8,
                     ymin=-4,
                     width = 0.35\textwidth,
                     height = 0.35\textwidth / \aspectratio,
                     scale only axis=true,
                     ymax=5,
                     ticks=none,
                     legend cell align={left}
                     ]

            \addplot [thick, black] table{\datadir/example_function.csv\ext};

            \addplot [thick, dashed, aux_color] table{\datadir/example_function_envelope.csv\ext};

            \legend{$f(x)$, $f^{**}(x)$}

        \end{axis}
    \end{tikzpicture}
\fi
\caption{A nonconvex function $f$ (solid black)
and its convex envelope $f^{**}$ (dashed \auxcolorname{}).}
\label{f-conv-env-example}
\end{figure}

The convex envelope can be defined several 
other equivalent ways.  For example,
$f^{**}$ is the greatest convex function that is a global
underestimator of $f$.
It is also the (Fenchel) conjugate of the conjugate of $f$, \ie,
$(f^*)^*$, where the superscript $*$ is the traditional notation
for the conjugate function.   (This explains why we denote
the convex envelope of $f$ as $f^{**}$.)
An example is shown in figure~\ref{f-conv-env-example}.

If $f(x)$ is convex when restricted to $x \leq 0$ and also
when restricted to $x \geq 0$, we can require that $v \geq 0$ and $w \leq 0$ 
in~(\ref{e-conv-env-def}),
\ie, we can define the convex envelope as
\begin{equation}\label{e-signed-conv-env-def}
f^{**}(x) = 
\inf \{\theta f(v) + (1-\theta) f(w) \mid 
\theta \in [0,1],~x = \theta v
+ (1 - \theta) w, ~v \geq 0,~w\leq 0\}.
\end{equation}

\subsection{Convex relaxation with borrowed curvature}
In this section we describe a convex relaxation of the TAM problem~(\ref{e-prob})
using the convex envelope.
We first eliminate the post-trade holdings variable $h$ to 
express the problem in terms of the trade list $u$,
with an objective that is the sum of a separable function and one
that is not separable.
This form is:
\begin{equation}\label{e-prob2}
\begin{array}{ll}
\text{maximize} &
-f_0(u) - \sum_{i=1}^n f_i(u_i)\\
\text{subject to} & u \in \tilde {\mathcal U} \\
\end{array}
\end{equation}
with variable $u\in \reals^n$, where the constraint set is
\[
\tilde {\mathcal U} =  
\{u \in \mathcal U \mid h_{\rm init} + u \in \mathcal H,
~\ones^Tu = c_{\rm init} - c_{\rm des}\}.
\]
The constraint set $\tilde {\mathcal U}$ includes the original constraint
$u\in \mathcal U$ as well as the holdings constraint $h\in \mathcal H$
and the post-trade cash constraint, and is convex.
The non-separable part of the objective function is
\[
f_0(u) = \gamma_{\rm risk} 
(h_{\rm init} - h_b + u)^T X \Sigma X^T (h_{\rm init} - h_b + u),
\]
which is the systematic component of risk~(\ref{e-factor-risk}), and is convex.  
The separable part corresponding to asset $i$ is
\begin{equation}\label{e-nonsep-sep-form}
f_i(u_i) = -\alpha_i u_i  
+ \gamma_{\rm risk} D_{ii} (h_{{\rm init}, i} - h_{b, i} + u_i)^2
+ \gamma_{\rm tc} \kappa_i |u_i|
+ \gamma_{\rm tax} L_i(u_i),
\end{equation}
which includes contributions from the expected return, 
specific risk~(\ref{e-idio-risk}), transaction cost, and tax liability.

The functions $f_i$ are
piecewise quadratic and nonconvex in general,
but are convex when $u_i \leq 0$ or $u_i \geq 0$. 
The problem~(\ref{e-prob2}), which is equivalent to the original
TAM problem~(\ref{e-prob}), is not convex because the functions $f_i$
are not convex.  However, if we fix the sign of each $u_i$,
the problem~(\ref{e-prob2}) becomes convex, and therefore easy to solve.
(In fact it suffices to fix the sign of $u_i$ for each asset where we hold 
at least one lot at a loss; the other $f_i$ are convex.)

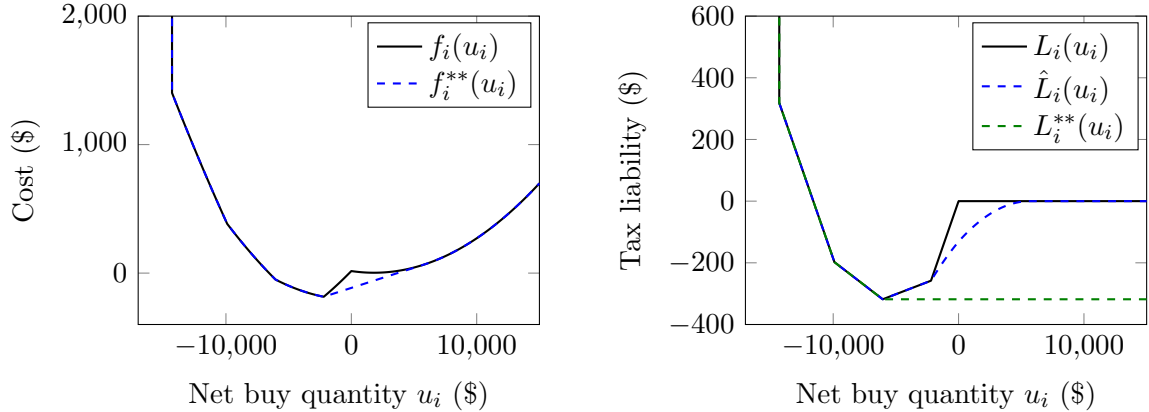
\begin{figure}
\centering
        \centering

        \ifmakeplots
            \begin{tikzpicture}
                \begin{axis}[xlabel = Net buy quantity $u_i$ (\$),
                             ylabel = Cost (\$),
                             xmin=-17000,
                             xmax=15000,
                             ymin=-400,
                             ymax=2000,
                             width = 0.35\textwidth,
                             height = 0.35\textwidth / \aspectratio,
                             scale only axis=true,
                             xticklabel style={/pgf/number format/fixed},
                             scaled ticks = false,
                             xtick={-10000, 0, 10000},
                             legend cell align={left}
                             ]

                    \addplot [thick, black] table{\datadir/separable_cost_function.csv\ext};

                    \addplot [thick, dashed, soph_heur_color] table{\datadir/separable_cost_function_envelope.csv\ext};

                    \legend{$f_i(u_i)$, $f_i^{**}(u_i)$}

                \end{axis}
            \end{tikzpicture}
        \fi
    %
    %
    \hfill
        \ifmakeplots
            \begin{tikzpicture}
                \begin{axis}[xlabel = Net buy quantity $u_i$ (\$),
                             ylabel = Tax liability (\$),
                             xmin=-17000,
                             xmax=15000,
                             ymin=-400,
                             ymax=600,
                             width = 0.35\textwidth,
                             height = 0.35\textwidth / \aspectratio,
                             scale only axis=true,
                             xticklabel style={/pgf/number format/fixed},
                             scaled ticks = false,
                             xtick={-10000, 0, 10000},
                             legend cell align={left}
                             ]

                    \addplot [thick, black] table{\datadir/tax_liability_function_with_loss.csv\ext};

                    \addplot [thick, dashed, soph_heur_color] table{\datadir/tax_liability_function_approximate.csv\ext};

                    \addplot [thick, dashed, basic_heur_color] table{\datadir/tax_liability_function_envelope.csv\ext};

                    \legend{$L_i(u_i)$, $\hat L_i(u_i)$, $L_i^{**}(u_i)$}

                \end{axis}
            \end{tikzpicture}
        \fi
\caption{
\emph{Left.} $f_i(u_i)$ (solid black) and its convex envelope $f_i^{**}(u_i)$ 
(\sophcolorname{} dashed).
\emph{Right.} Nonconvex tax liability function $L_i(u_i)$ (solid black),
its convex envelope $L_i^{**}(u_i)$ (dashed \basiccolorname{} line), and the 
approximation used in our sophisticated relaxation $\hat L_i(u_i)$ (dashed \sophcolorname{}).
Even though $\hat L_i(u_i)$ is nonconvex, we can still solve the problem
globally and efficiently.  }
\label{f-approx-tax-function}
\label{f-separable-function-envelope}
\end{figure}

\paragraph{Relaxed TAM problem with borrowed curvature.}
We can now obtain a convex relaxation of the TAM problem
by replacing $f_i$ with $f_i^{**}$, to obtain
\begin{equation}\label{e-prob-relax2}
\begin{array}{ll}
\text{maximize} & 
-f_0(u) - \sum_{i=1}^n f_i^{**}(u_i)\\
\text{subject to} &
u \in \tilde{\mathcal U}.
\end{array}
\end{equation}
This is a convex problem, which can be formulated as a second-order cone problem (SOCP)
\citep[\S4.4.2]{cvxbook}.
The convex envelopes $f_i^{**}$ are convex and also piecewise quadratic.
Figure~\ref{f-separable-function-envelope} plots $f_i$ and $f_i^{**}$.


The objective of the relaxation~(\ref{e-prob-relax2}) is an upper bound
on the objective of the original TAM problem.  
It follows that its optimal
objective value $U^\star_{\rm relax}$ is an upper bound on the
optimal value of the original TAM problem.
It follows that its optimal
objective value $U^\star_{\rm relax}$ is an upper bound on the
optimal value of the original TAM problem, \ie,
\[
U^\star \le U_{\rm relax}^\star.
\]

The gap $U^\star_{\rm relax} - U^\star$ can be bounded in terms of $k$,
the number of factors in the risk model,
and the distances between the separable functions $f_i$ and their convex envelopes $f_i^{**}$.
This is an application of the Shapley--Folkman Lemma \citep{bertsekas1982,udell2016bounding}.

\paragraph{TAM problem with approximate tax liability.}
By re-introducing the post-trade holdings variable,
the relaxed problem~(\ref{e-prob-relax2}) can be written as
\begin{equation}
\label{e-approx-tam-prob}
\begin{array}{ll}
\text{maximize} &
\alpha^T u
- \gamma_\text{risk} (h - h_b)^T V (h - h_b)
- \gamma_\text{tc} \kappa^T |u|
- \gamma_\text{tax} \hat L(u)  \\
\text{subject to} &
h = h_{\rm init} + u, \quad
\ones^T u =c_\text{init}-c_\text{des} \\
& u \in \mathcal U, \quad h \in \mathcal H.
\end{array}
\end{equation}
with decision variables $u$ and $h$.
This is the TAM problem with
the tax liability functions $L$ replaced by \emph{approximate} tax liability function $\hat L$,
defined as
\[
\hat L(u) = \sum_{i=1}^n \hat L_i(u_i),
\]
where $\hat L_i = L_i + (f_i^{**} - f_i)/\gamma_{\rm tax}$.
Problem~(\ref{e-approx-tam-prob}) can be solved exactly using convex optimization,
even though the functions $\hat L_i$ are not convex.
This is possible because the nonconvex function $\hat L$
borrows curvature from the other separable objective terms,
resulting in an objective function that is concave.
(See figure~\ref{f-separable-function-envelope}.)

\section{Approximate solution methods}
\label{s-approx-sol}

This section describes heuristic solution methods for the TAM problem.
The methods involve solving two convex 
optimization problems in two stages.
\begin{enumerate}
    \item \emph{Guess the vector $z$ of signs of an optimal $u$.}
        This is done by solving the relaxation~\eqref{e-prob-relax2} of the TAM problem (which 
        in addition provides an upper bound on $U^\star$).
    \item \emph{Solve TAM with these sign constraints.}
        Add the sign constraints $z_iu_i \geq 0$, $i=1, \ldots, n$ to
        the TAM problem and solve.
        With these constraints,
        the TAM problem is a convex QP and can be efficiently solved.
\end{enumerate}
There are several choices for step~1, which we describe below.
We only need to specify the sign of $u_i$ for 
assets in which we hold at least one lot at a loss.

\paragraph{Methods for guessing the sign.}
For step~1, we solve or the relaxation~(\ref{e-prob-relax2}).
There are several choices for guessing the signs of $u_i$ from the 
solution of one of this relaxation.
The most obvious method is to use the sign of 
the solution of the relaxation, \ie, $z = \sign(u^\star_{\rm relax})$.

A less obvious method is a random choice of the signs with probabilities
taken from the solution of the relaxation.
For each $i$, we obtain the values $\theta_i$ in the 
convex envelope definition~(\ref{e-signed-conv-env-def}) 
for $f_i$.
We then set $z_i = 1$ with probability $\theta_i$,
and $z_i = -1$ with probability $1 -\theta_i$.
(This is done independently for each $i$.)
Thus we use the values in the convex envelope as probabilities 
on whether we buy or sell each asset.
This method can be used to generate multiple candidate sign vectors,
and we can compare the objectives after step~2
and use the one with the largest objective.

In many numerical experiments we found that the method that
performs best is to solve the relaxation~(\ref{e-prob-relax2}),
and then use the randomized method to guess a set of signs.  
(This is despite the fact that the simple rounding method
is guaranteed to produce feasible sign constraints for the TAM problem,
and the randomized method is not.)
We have also found that generating multiple sets of candidate signs does not 
substantially improve the results.  This method requires two convex
optimization solves: one to solve the relaxation (an SOCP), and one to solve
the orignal TAM problem with the sign of $u$ fixed (a QP).

\section{Numerical examples}\label{s-examples}
We demonstrate these methods by simulating a tax-loss harvesting strategy,
in which we solve the TAM problem once a month to generate the trade list.
First, we show a six-year backtest of such a strategy.
Then, we use this backtest (and others like it)
to generate realistic instances of the TAM problem,
which we use to evaluate the methods of section~\ref{s-approx-sol}.

\subsection{Benchmark and data}

All of our simulations use the S\&P 500 as the benchmark, with data over
the period 2002 to 2019.
Our universe includes all assets that were in the S\&P 500 at any point
over that time interval, which gives $n=998$.
We included a constraint that we only purchase shares of current S\&P 500 constituents.
This prevents us from purchasing assets that, at the time of the simulated trade, have never been in the benchmark.
It also means we don't increase our holdings of former S\&P 500 constituents
(but we also do not require them to be immediately sold).
If any asset is delisted,
we liquidate the asset immediately, incurring the associated tax liability.

We take $\alpha=0$, \ie, we do not have any views on the active returns,
so our goal is to simply track the benchmark portfolio while minimizing
tax liability.
Our risk model parameters $\Sigma$, $X$, and $D$
are from the Barra US Equity model \citep{menchero2011barra},
which uses $k=72$ factors.
Our cash target $c_{\rm des}$ is given by~(\ref{e-cdes-choice}) with $\eta =0.005$,
\ie, we hold 50 basis points in cash after each trade.
We use tax rates $\rho_{\rm lt} = 0.238$ and $\rho_{\rm st} = 0.408$,
which reflect the current highest marginal tax rates in the United States
for long-term and short-term capital gains, respectively.
We used the conservative value $\kappa_i = 0.0005$ for all transaction
costs, \ie, the bid-ask spread is 10 basis points for all assets.
The parameter $\gamma_{\rm risk}$ was scaled with the account value,
so that $\gamma_{\rm risk} = \tilde \gamma_{\rm risk} (\ones^T h_{\rm init} + c_{\rm init})$,
with $\tilde \gamma_{\rm risk} = 200$.
The other trade-off parameters were $\gamma_{\rm tc} = 1$ and $\gamma_{\rm tax} = 1$.

\subsection{Backtests}
Our dataset consists of 204 months over a 17 year period from August 2002 through August 2019.
We use this dataset to carry out 12 different, staggered six-year-long backtests.  The first one
starts in August 2002 and ends in July 2008; the last one starts in August 2013 and 
ends in July 2019.
In these backtests, monthly trading means we trade on the first business day
more than 31 calendar days after the last trade.
For each trade, the initial cash amount $c_{\rm init}$
is adjusted for the realized transaction cost $\kappa^T |u|$ of the last trade,
as well as cash inflows due to dividends and other corporate actions.
In the backtests, we round the trade lists to an integer number of shares.
Each backtest starts with a portfolio of \$1M in cash.

Each month, the trade list is determined by solving the TAM problem using one of two methods:
\begin{itemize}
\item \emph{Heuristic.} We solve the relaxation~(\ref{e-prob-relax2})
    and use the randomized rounding method.
\item \emph{Mixed-integer method.} We use the mixed-integer mode of CPLEX (version 12.9) to solve the TAM problem directly,
    with a time limit of 300 seconds.
\end{itemize}
For the heuristic, we used CPLEX (as a QP/SOCP solver) to solve
the convex relaxations and to generate the final trade list.

\begin{figure}
\centering
\ifmakeplots
    \begin{tikzpicture}

        \begin{groupplot}[group style={group size=1 by 2},
                          width = 0.85\textwidth,
                          height = 0.3\textwidth,
                          date coordinates in = x,
                          xmin = 2012-08-15,
                          xmax = 2018-08-15,
                          legend cell align = {left},
                          scaled ticks = false,
                          xtick = {2013-01-01, 2014-01-01, 2015-01-01, 2016-01-01, 2017-01-01, 2018-01-01}, 
                          yticklabel style={/pgf/number format/fixed},
                          xticklabel=\year,
                          clip mode = individual ] 

        \nextgroupplot[ylabel = Active risk (\%),
                       ymin = 0.00,
                       ymax = 0.60]

            \addplot [thick, black, on layer=foreground] table{\datadir/backtest_risk_mip.csv\ext};

            \addplot [thick, soph_heur_color, on layer=main] table{\datadir/backtest_risk_heuristic.csv\ext};

        \nextgroupplot[ylabel = Cumulative tax liability (\$),
                       xlabel = Date,
                       ymin = -60000,
                       ymax =  0]

            \addplot [thick, black, on layer=foreground] table{\datadir/backtest_tax_mip.csv\ext};
            \addlegendentry{Mixed integer method}

            \addplot [thick, soph_heur_color, on layer=main] table{\datadir/backtest_tax_heuristic.csv\ext};
            \addlegendentry{Heuristic}

        \end{groupplot}
    \end{tikzpicture}
\fi
\caption{
The active risk (top) and cumulative tax liability (bottom) of a backtest for both solution methods.
}
\label{f-backtest-risk}
\end{figure}
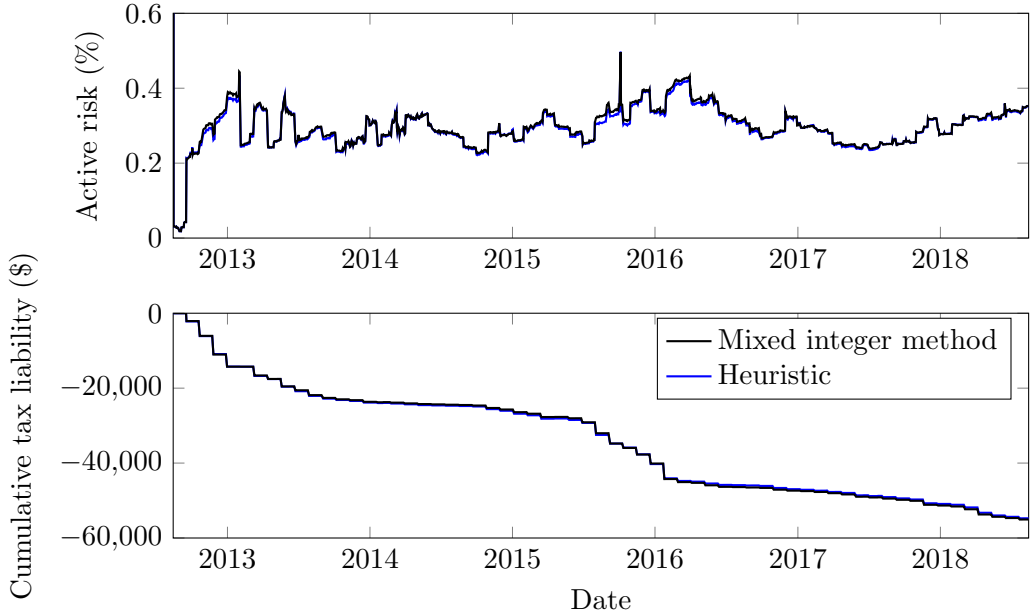

\paragraph{Example backtest.}
Figure~\ref{f-backtest-risk} shows the results of one of our backtests, initiated in August 2013.
The top plot shows the active risk, and the bottom plot shows the 
cumulative tax liability,
which is the net realized gain, accounting for long- and short-term tax rates.
(This quantity is negative, meaning we are realizing a net loss).
Here we use the conventional definition of active risk,
which is the square root of the definition given in section~\ref{s-TAM}
and is scaled down by the account value.

These results show that a tax-aware trading scheme can indeed 
track a benchmark while simultaneously realizing capital losses.
It is interesting to note that losses are harvested even during bull markets.
The rate of tax-loss harvesting decreases with the life of the fund,
since more lots are held at a gain. 
We note that the two methods have nearly indistinguishable performance in the backtest.
The backtests we ran over other time windows had similar outcomes.

\subsection{Detailed comparison of solution methods}
This section compares the performance of both solution methods.
We use the data from 12 backtests, each six years long,
giving a total of 744 instances of the TAM problem.
(We exclude the initial trade, in which the account holds only cash.)
For these problem instances,
we compute the utility achieved by the heuristic method, denoted $U_{\rm relax,round}$,
the utility of the mixed integer method, denoted $U_{\rm mip}$,
as well as the upper bound $U_{\rm relax}^\star$.
To make the utility~(\ref{e-utility}) comparable across the problem instances,
we divide it by the account value $\ones^T h_{\rm init} + c_{\rm init}$.
(This number is the monthly after-tax expected return
adjusted for risk and transaction costs
and is measured in percent or basis points.)
The mean optimal utility $U^\star$ across these 744 problems
is $2$ basis points, and the standard deviation is $55$ basis points.

\paragraph{Evaluation of the heuristic.}
The results are shown in table~\ref{t-results}.
We see that for 678 of the 744 instances, 
$U_{\rm relax,round} = U_{\rm relax}^\star$,
\ie, the heuristic solves the TAM problem and provides a certificate of optimality,
to within the solver numerical accuracy, which is around $0.05$ bp.
We also see that the (mean) differences in utility obtained by the heuristic and the bound are very close,
\ie, within fractions of a basis point.
To summarize, the heuristic
produces an optimal trade list (and certificate of optimality)
for the vast majority of the problem instances.
We note that for all 744 problem instances,
the heuristic is never suboptimal by more than a few basis points.

We now compare the heuristic to the mixed-integer method.
Note that if CPLEX solves the mixed-integer problem within the $300$ second time limit,
we have $U_{\rm mip} = U^\star$.
The mixed-integer method times out (and therefore, does not necessarily globally solve the problem)
in 179 of the 744 cases.
Among the 565 instances that the mixed-integer method solves (within 300 seconds),
in 549 instances the 
heuristic method also solves the problem to within $0.05$ basis points,
and is never more than $0.3$ basis points suboptimal in the remaining 16 cases.
For the 179 cases in which the mixed integer method fails to solve the problem,
we observe that the 
heuristic achieves a normalized utility within $0.05$ basis points
of the mixed integer method in 89 cases,
and outperforms it by more than $0.05$ basis points in 68 cases (by up to $2$ basis points).
In only 22 of the 179 instances did the 
heuristic underperform the mixed integer method 
by more than $0.05$ basis points,
and never underperformed by more than $2$ basis points.

\begin{table}
    \footnotesize
    \centering
    \begin{tabular}{lccdd}
        \toprule
        & \multicolumn{1}{c}{Improved vs MIP} & \multicolumn{1}{c}{Solved}  
        & \multicolumn{1}{c}{Mean subopt.\ gap} & \multicolumn{1}{c}{Mean MIP gap} \\
        & \multicolumn{1}{c}{$U\ge U_{\rm mip}$} & \multicolumn{1}{c}{$U=U^\star_{\rm relax}$}  
        & \multicolumn{1}{c}{$U^\star_{\rm relax} - U$} & \multicolumn{1}{c}{$U^\star_{\rm mip} - U$} \\
        & \multicolumn{1}{c}{(count)} & \multicolumn{1}{c}{(count)}  
        & \multicolumn{1}{c}{(bp)} & \multicolumn{1}{c}{(bp)} \\
        \midrule
        Heuristic     & 706 & 678 & 0.02     & -0.01 \\
        MIP (300s)        & 744 & 646 & 0.03     & 0   \\
        \bottomrule
    \end{tabular}
    \caption{
    Comparison of the solution methods.
    The first two columns show the number of problem instances (out of 744 total)
    for which the achieved utility $U$ of each method exceeds the mixed-integer method's utility $U_{\rm mip}$
    (first column)
    or matches the upper bound $U^\star_{\rm relax}$ exactly (second column).
    Equality here is to within $0.05$ bp, the approximate solver tolerance.
    The last two columns give the average difference between the achieved utility $U$
    and $U_{\rm mip}$ (third column) or $U^\star_{\rm relax}$ (fourth column).
    }
    \label{t-results}
\end{table}

\begin{figure}
    \centering
    \ifmakeplots
        \begin{tikzpicture}[/pgfplots/scale only axis,
                /pgfplots/width=6cm,
                /pgfplots/height=6cm
            ]

            \def\timeout{300}
            \def\mintime{0.1}
            \def\maxtime{500}
            \def\logmintime{-1}
            \def\logmaxtime{2.7}
            \ifjota
                \def\histwidth{.4\columnwidth}
                \def\sidehistheight{.15\columnwidth}
            \else
                \def\histwidth{.5\columnwidth}
                \def\sidehistheight{.2\columnwidth}
            \fi

            \begin{axis}[
                 name = main axis,
                 xmin = \mintime,
                 xmax = \maxtime,
                 ymin = \mintime,
                 ymax = \maxtime,
                 xmode = log,
                 ymode = log,
                 width = \histwidth,
                 height = \histwidth,
                 xlabel = {Solve time, mixed integer (s)},
                 ylabel = {Solve time, heuristic (s)},
                 xticklabel style={/pgf/number format/fixed},
                 yticklabel style={/pgf/number format/fixed},
                 scaled ticks = false,
                 xtick = {0.1, 1, 10, 100, 1000},
                 ytick = {0.1, 1, 10, 100, 1000},
                 xticklabels = {0.1, 1, 10, 100, },
                 yticklabels = {   , 1, 10, 100, },
                 log ticks with fixed point
            ]
            \addplot [only marks, mark size=0.5] table[x index=1, y index=2] {\datadir/solve_times.csv\ext};
            \addplot [dashed, black] coordinates {(\mintime, \mintime)  (\maxtime, \maxtime)};
            \addplot [dashed, red] coordinates {(\timeout, \mintime)  (\timeout, \maxtime)};
            \end{axis}

            \begin{axis}[
                anchor = south west,
                at = (main axis.north west),
                height = \sidehistheight,
                width = \histwidth,
                xtick = \empty,
                ytick = \empty,
                ymin = 0,
                xmin = \logmintime,
                xmax = \logmaxtime
            ]
            \addplot [
                hist={data=x, bins=30, data min = \logmintime, data max = \logmaxtime},
                fill=gray!50
            ] table[x index=1] {\datadir/solve_times_log.csv\ext};
            \end{axis}

            \begin{axis}[
                anchor = north west,
                at = (main axis.north east),
                width = \sidehistheight,
                height = \histwidth,
                xtick = \empty,
                ytick = \empty,
                xmin = 0,
                ymin = \logmintime,
                ymax = \logmaxtime
            ]
            \addplot [
                hist = {handler/.style={xbar interval}, bins=30, data min=\logmintime, data max=\logmaxtime},
                x filter/.code=\pgfmathparse{rawy}, 
                y filter/.code=\pgfmathparse{rawx}, 
                fill = gray!50,
            ] table[y index=2] {\datadir/solve_times_log.csv\ext};
            \end{axis}
        \end{tikzpicture}
    \fi
    \caption{
        The algorithm run times of the 720 problem instances using the relax-and-round heuristic
        and the mixed-integer solution, with each problem instance shown as a single dot.
        The dashed red line shows the maximum allowed time of the mixed-integer solver.
    }
    \label{f-solve-time-scatterplot}
\end{figure}
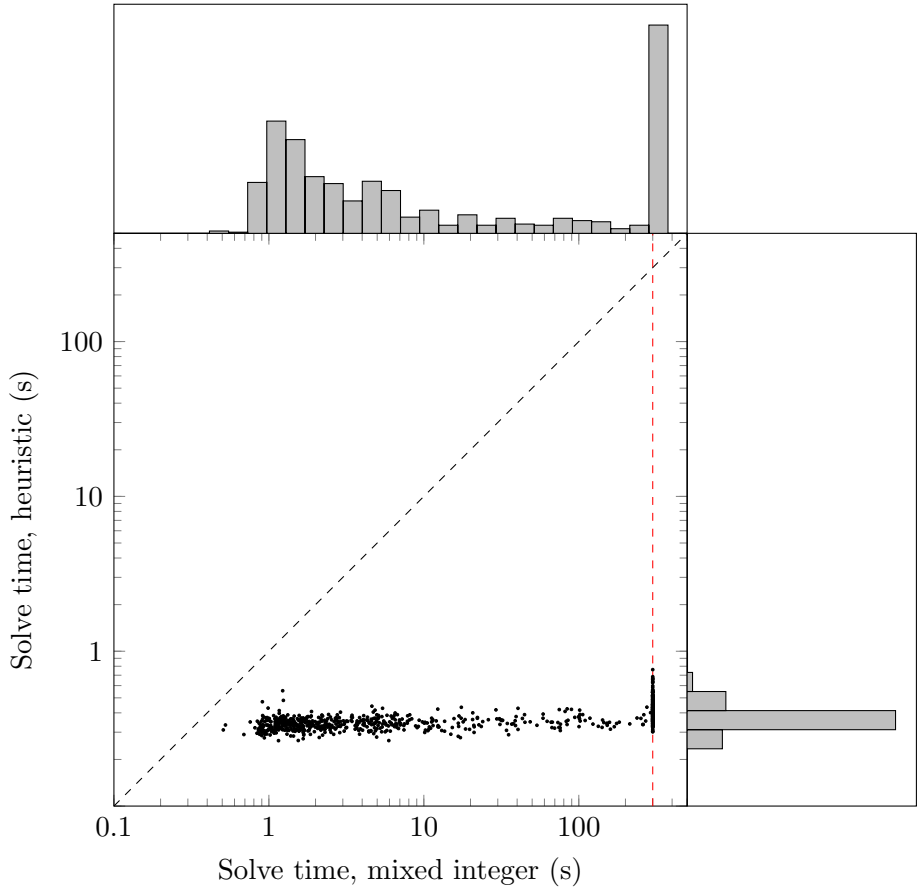

\paragraph{Algorithm run times.}
Figure~\ref{f-solve-time-scatterplot}
shows the algorithm run times for the 
heuristic and the mixed integer method on a scatter plot.
All of the points are below the dashed black line,
which indicates that the heuristic method was faster in all cases.
Out of the 744 problem instances, 179 took $300$ seconds using the mixed-integer method,
which was the maximum time allowed.

\section{Conclusion}
We formulate tax-aware portfolio construction as a nonconvex optimization problem,
and we present a heuristic for this problem based on convex optimization.
This method is reliably fast:
for problems with several hundred assets and several dozen factors,
it takes less than a second.
We compare our heuristic against
the standard, mixed-integer quadratic programming formulation, solved using CPLEX,
on realistic problem instances.
When the mixed-integer method is limited to five minute run times,
we find that our heuristic outperforms it more often than not,
despite being several hundred times faster.
This speed is not necessary for monthly (or even daily) trading,
but is useful for backtesting and Monte Carlo simulation,
possibly over hundreds of thousands of individualized accounts.
Our method also produces a bound on the optimal value.
For realistic data, the bound usually tight enough
that it certifies that the heuristic solved the problem globally.
In future work, we will extend this method to other nonconvex terms that are often present 
in practical portfolio optimization problems.

\paragraph{Acknowledgements.}
We would like to thank Emmanuel Cand\`{e}s for useful discussions and feedback.
We would also like to thank Eric Kisslinger for identifying an important error 
in an early version of the software.

\printbibliography

\ifjota
\else

    \appendix

    \section{SOCP formulation}\label{s-socp-formulation}
    Here we explain how to represent the convex envelope $f_i^{**}$ in a cone program,
    by expressing its epigraph using a cone representation, as described by \citet{gb08}.
    The technique given here are similar to those used to represent perspectives
    of convex functions \citep[\S~2]{perspective}.

    Consider a function $f:\reals\to \reals \cup \{\infty\}$, of the form
    \[
    f(x) = \left\{ \begin{array}{ll} f_-(x) & x < 0\\ f_+(x) & x \geq 0,
    \end{array}\right.
    \]
    where $f_-$ and $f_+$ are both convex, with $f_-(x)=+\infty$ for $x>0$
    and $f_+(x)=+\infty$ for $x<0$.
    We assume that each of these functions has a so-called cone representation.
    This means that $f_-(x)$ is the optimal value of a cone program
    \[
    \begin{array}{ll} 
    \mbox{minimize} & c_-^T z_- \\
    \mbox{subject to} & A_-(x,z_-) = b_-, \quad (x,z_-) \in \mathcal K_-,
    \end{array}
    \]
    with variable $z_-$, where $\mathcal K_-$ is a cone.  We assume a similar 
    representation for $f_+$.

    Our goal is to represent the convex envelope~(\ref{e-signed-conv-env-def})
    as the optimal value of a cone program.
    Using the cone representations of $f_-$ and $f_+$, we can express $f^{**}(x)$
    as the optimal value of the problem
    \[
    \begin{array}{ll} 
    \mbox{minimize} & \theta c_-^T z_- + (1-\theta) c_+^T z_+ \\
    \mbox{subject to} & A_-(v,z_-) = b_-, \quad (v,z_-) \in \mathcal K_-,\\
    & A_+(w,z_+) = b_+, \quad (w,z_+) \in \mathcal K_+,\\
    & x= \theta v + (1-\theta) w,\\
    & v \geq 0, \quad w \leq 0, \quad 0 \leq \theta \leq 1,
    \end{array}
    \]
    with variables $\theta$, $z_-$, $z_+$, $v$, and $w$.
    The objective terms and the equality constraint involving $x$ 
    contain the product of two variables, and is not convex.

    We will now change variables to obtain an equivalent convex problem.
    Define the variables
    \begin{align}
    \label{e-change-of-variables}
    \tilde z_- = \theta z_-, \quad
    \tilde v = \theta v, \quad
    \tilde z_+ =(1-\theta) z_+, \quad
    \tilde w = (1-\theta)w.
    \end{align}
    We can express the problem above using these variables, and the original 
    variable $\theta$, as
    \begin{align}
    \label{e-equivalent-socp}
    \begin{array}{ll} 
    \mbox{minimize} & c_-^T \tilde z_- + c_+^T \tilde z_+ \\
    \mbox{subject to} & A_-(\tilde v,\tilde z_-) = \theta b_-, \quad 
    (\tilde v,\tilde z_-) \in \mathcal K_-,\\
    & A_+(\tilde w,\tilde z_+) = (1-\theta) b_+, \quad (\tilde w,\tilde z_+)
    \in \mathcal K_+,\\
    & x= \tilde v + \tilde w,\\
    & \tilde v \geq 0, \quad \tilde w \leq 0, \quad 0 \leq \theta \leq 1,
    \end{array}
    \end{align}
    with variables $\theta$, $\tilde z_-$, $\tilde z_+$, $\tilde v$, and $\tilde w$.
    This problem is jointly convex in all variables, and $x$, so it is a cone 
    representation of $f^{**}$.

    We note that for the change of variables~(\ref{e-change-of-variables}) to be invertible, 
    we must include in problem~(\ref{e-equivalent-socp})
    the constraint that $\tilde z_-$ must be $0$ if $\theta$ is $0$.
    Because this additional constraint only restricts points on the boundary
    of the feasible set of problem~(\ref{e-equivalent-socp}),
    we can safely ignore it without changing the optimal value of the problem,
    assuming Slater's condition holds.
    Similar arguments apply for $\tilde v$, $\tilde z_+$, and $\tilde w$.

    For the specific case where $f$ is piecewise quadratic,
    (\eg, the separable cost functions $f_i$ given in (\ref{e-nonsep-sep-form})),
    the cone representations of $f_-$, $f_+$, and $f^{**}$ are second-order cone programs (SOCPs).
    This means the relaxation~(\ref{e-prob-relax2}) can be expressed as an SOCP.
\fi

\end{document}